\documentclass[12pt]{article}

\usepackage{amsmath,amssymb,amsbsy,amsfonts,amsthm,
            latexsym,amsopn,amstext,amsxtra,euscript,amscd,mathrsfs}

\begin{document}

\newtheorem{theorem}{Theorem}
\newtheorem{lemma}{Lemma}

\newtheorem{proposition}[theorem]{Proposition}
\newtheorem{corollary}[theorem]{Corollary}

\theoremstyle{definition}
\newtheorem*{definition}{Definition}
\newtheorem*{remark}{Remark}
\newtheorem*{example}{Example}


\def\cA{\mathcal A}
\def\cB{\mathcal B}
\def\cC{\mathcal C}
\def\cD{\mathcal D}
\def\cE{\mathcal E}
\def\cF{\mathcal F}
\def\cG{\mathcal G}
\def\cH{\mathcal H}
\def\cI{\mathcal I}
\def\cJ{\mathcal J}
\def\cK{\mathcal K}
\def\cL{\mathcal L}
\def\cM{\mathcal M}
\def\cN{\mathcal N}
\def\cO{\mathcal O}
\def\cP{\mathcal P}
\def\cQ{\mathcal Q}
\def\cR{\mathcal R}
\def\cS{\mathcal S}
\def\cU{\mathcal U}
\def\cT{\mathcal T}
\def\cV{\mathcal V}
\def\cW{\mathcal W}
\def\cX{\mathcal X}

\def\Z{{\mathbb Z}}
\def\R{{\mathbb R}}
\def\F{{\mathbb F}}
\def\N{{\mathbb N}}
\def\C{{\mathbb C}}
\def\Q{{\mathbb Q}}
\def\P{{\mathbb P}}

\def\fc{\mathcal F} 

\def\Fp{\F_p}
\def\e{\mathbf{e}}
\def\ep{\mathbf{e}_p}
\def\eps{{\varepsilon}}
\def\mand{\qquad \mbox{and} \qquad}
\def\scr{\scriptstyle}
\def\\{\cr}
\def\({\left(}
\def\){\right)}
\def\[{\left[}
\def\]{\right]}
\def\<{\langle}
\def\>{\rangle}
\def\fl#1{\left\lfloor#1\right\rfloor}
\def\rf#1{\left\lceil#1\right\rceil}
\def\le{\leqslant}
\def\ge{\geqslant}
\def\sL{\mathscr L}
\def\sM{\mathcal M}
\def\sT{\mathcal T}
\def\ds{\displaystyle}
\renewcommand{\mod}{\textrm{mod}\,}

\def\xxx{\vskip5pt\hrule\vskip5pt}

\newcommand{\comm}[1]{\marginpar{%
\vskip-\baselineskip 
\raggedright\footnotesize
\itshape\hrule\smallskip#1\par\smallskip\hrule}}

\title{\sc On primitive Dirichlet characters and the Riemann hypothesis}

\author{
{\sc William D.~Banks} \\
{Department of Mathematics} \\
{University of Missouri} \\
{Columbia, MO 65211 USA} \\
{\tt bbanks@math.missouri.edu} \\
\and
{\sc Ahmet M.~G\"ulo\u glu} \\
{Department of Mathematics} \\
{University of Missouri} \\
{Columbia, MO 65211 USA} \\
{\tt ahmet@math.missouri.edu} \\
\and
{\sc C.~Wesley Nevans} \\
{Department of Mathematics} \\
{University of Missouri} \\
{Columbia, MO 65211 USA} \\
{\tt nevans@math.missouri.edu}}

\date{}

\maketitle

\begin{abstract}
For any natural number $n$, let $\cX'_n$ be the set of primitive
Dirichlet characters modulo~$n$. We show that if the Riemann
hypothesis is true, then the inequality $|\cX'_{2n_k}|\le
C_2\,e^{-\gamma}\,\varphi(2n_k)/\log\log(2n_k)$ holds for all
$k\ge 1$, where $n_k$ is the product of the first $k$ primes,
$\gamma$ is the Euler-Mascheroni constant, $C_2$ is the twin prime
constant, and $\varphi(n)$ is the Euler function.  On the other
hand, if the Riemann hypothesis is false, then there are
infinitely many $k$ for which the same inequality holds and
infinitely many $k$ for which it fails to hold.
\end{abstract}

\section{Introduction}
\label{sec:intro}

For any natural number $n$, let $\cX_n$ be the set of Dirichlet
characters modulo~$n$, and let $\cX'_n$ be the subset of
\textit{primitive} characters in $\cX_n$.

The purpose of the present note is to establish a connection
between the \textit{classical} Riemann hypothesis and the
collection of sets $\{\cX'_n:n\in\N\}$.  Our work is motivated by
and relies on the 1983 paper of J.-L.~Nicolas~\cite{Nic} in which
a relation is established between the Riemann hypothesis and
certain values of the Euler function $\varphi(n)$; see
also~\cite{Robin}.

\begin{theorem}
\label{thm:main} For every $k\ge 1$, let $n_k$ be the product of
the first $k$ primes. Let $\gamma$ be the Euler-Mascheroni
constant and $C_2$ the twin prime constant.

\begin{itemize}

\item[$(i)$] If the Riemann hypothesis is true, then the
inequality
\begin{equation}
\label{eq:RHtrue} |\cX'_{2n_k}|\le
C_2\,e^{-\gamma}\,\frac{\varphi(2n_k)}{\log\log(2n_k)}
\end{equation}
holds for all $k\ge 1$.

\item[$(ii)$] If the Riemann hypothesis is false, then there are
infinitely many $k$ for which \eqref{eq:RHtrue} holds and
infinitely many $k$ for which it fails to hold.
\end{itemize}
\end{theorem}

We recall that
$$
\gamma=\lim_{n\to\infty}\(\,\sum_{m=1}^n\frac{1}{m}-\log
n\)=0.5772156649\cdots,
$$
and
$$
C_2=\prod_{p > 2} \frac{p(p-2)}{(p-1)^2}=0.6601618158 \cdots.
$$

To prove the theorem, we study the ratios
$$
\rho(n)=\frac{|\cX'_n|}{|\cX_n|}\qquad(n\in\N).
$$
Note that $\rho(n)$ is the \emph{proportion} of Dirichlet
characters modulo~$n$ that are primitive characters.  Since
$\rho(n)\le 1$ for all $n\in\N$, and $\rho(p)=1-1/(p-1)$ for every
prime $p$, it is clear that
$$
\limsup_{n\to\infty}\rho(n)=1.
$$
As for the minimal order, we shall prove the following:
\begin{equation}
\label{eq:liminf-n} \liminf_{\substack{n\to\infty\\n\not\equiv
2\;(\mod 4)}} \rho(n)\log\log n=C_2\,e^{-\gamma}.
\end{equation}
Note that natural numbers $n\equiv 2\pmod 4$ are excluded since
$\rho(n)=0$ for those numbers; see~\eqref{eq:rho-defn} below.

In Section~\ref{sec:two} we show that the inequalities
\begin{equation}
\label{eq:rho-ineq} \rho(2n_k)\log\log(2n_k)\le\rho(n)\log\log
n\qquad(n\not\equiv 2\;(\mod 4),~\omega(n)=k)
\end{equation}
hold for every fixed $k>1$, where $\omega(n)$ is the number of
distinct prime divisors of $n$, and we also show that
\begin{equation}
\label{eq:lim-2nk}
\lim_{k\to\infty}\rho(2n_k)\log\log(2n_k)=C_2\,e^{-\gamma}.
\end{equation}
Clearly, \eqref{eq:liminf-n} is an immediate consequence
of~\eqref{eq:rho-ineq} and~\eqref{eq:lim-2nk}.

Since $|\cX_n|=\varphi(n)$ for all $n\in\N$, the
inequality~\eqref{eq:RHtrue} is clearly equivalent to
\begin{equation}
\label{eq:RHtrueB} \rho(2n_k)\log\log(2n_k)\le C_2\,e^{-\gamma}.
\end{equation}
In Section~\ref{sec:three} we study this inequality using
techniques and results from~\cite{Nic}, and these investigations
lead to the statement of Theorem~\ref{thm:main}.

\bigskip

\noindent{\bf Acknowledgement.} The authors wish to thank Pieter
Moree for his careful reading of the manuscript and for several
useful comments.

\section{Small values of $\rho(n)$}
\label{sec:two}

The cardinality of $\cX_n$ is $\varphi(n)$, and that of $\cX'_n$
is
$$
|\cX'_n|=n\prod_{p\,\|\,n}\(1-\frac{2}{p}\)\prod_{p^2\,\mid\,n}
\(1-\frac{1}{p}\)^2
$$
(see, for example, \cite[\S\,9.1]{MontVau}); hence, it follows
that
\begin{equation}
\label{eq:rho-defn} \rho(n)=\frac{\varphi(n)}{n}\prod_{p\,\|\,n}
\frac{p(p-2)}{(p-1)^2}\qquad(n\in\N).
\end{equation}

Turning to the proof of~\eqref{eq:rho-ineq}, let $k>1$ be fixed,
and denote by $\cS$ the set of integers $n\not\equiv 2\pmod 4$
with $\omega(n)=k$. Let $p_1,p_2,\ldots$ be the sequence of
consecutive prime numbers. For each integer $j\in\{0,\ldots,k\}$,
let $\cS_j$ be the set of numbers $n\in\cS$ that have precisely
$j$ distinct prime divisors larger than~$p_k$. Since $\cS$ is the
union of the sets $\{\cS_j\}$, to prove~\eqref{eq:rho-ineq} it
suffices to show that the inequalities
\begin{equation}
\label{eq:Sjbound} \rho(2n_k)\log\log(2n_k)\le\rho(n)\log\log n
\qquad(n\in\cS_j)
\end{equation}
hold for every fixed $j\in\{0,\ldots,k\}$.

For any $n\in\cS_0$ we can write $n=2p_1^{\alpha_1}\cdots
p_k^{\alpha_k}$ with each $\alpha_j\ge 1$.
Using~\eqref{eq:rho-defn} and the fact that $2n_k=2p_1\cdots p_k$
we have
$$
\rho(2n_k)=\rho(n) \prod_{\substack{j=2\\ (\alpha_j\ge 2)}}^k
\frac{p_j(p_j-2)}{(p_j-1)^2}\le\rho(n).
$$
Since $2n_k\le n$ we also have $\log\log(2n_k)\le\log\log n$,
and~\eqref{eq:Sjbound} follows for $j=0$.

Proceeding by induction, let us suppose that~\eqref{eq:Sjbound}
has been established for some $j\in\{0,\ldots,k-1\}$. If $n'$ is
an arbitrary element of $\cS_{j+1}$, then $q\mid n'$ for some
prime $q>p_k$; note that $q\ge 5$ since $k>1$. Writing
$n'=q^\alpha m$ with $q\nmid m$, we have $\omega(m)=k-1$, hence
for at least one index $i\in\{1,\ldots,k\}$ the prime $p_i$ does
not divide $m$. Put $n=p_i^\beta m$, where $\beta=2$ if $p_i=2$
and $\beta=1$ otherwise.  Clearly, $n\in\cS_j$. Also, $n\le n'$
since $q>\max\{p_i,2^2\}$, and thus $\log\log n\le\log\log n'$. Finally,
using~\eqref{eq:rho-defn} we see that
$$
\frac{\rho(n')}{\rho(m)}=\left\{%
\begin{array}{ll}
    1-1/(q-1) & \quad\hbox{if $\alpha=1$,} \\
    1-1/q & \quad\hbox{if $\alpha\ge 2$,} \\
\end{array}%
\right.
$$
and
$$
\frac{\rho(n)}{\rho(m)}=\left\{%
\begin{array}{ll}
    1-1/(p_i-1) & \quad\hbox{if $\beta=1$,} \\
    1/2 & \quad\hbox{if $\beta=2$.} \\
\end{array}%
\right.
$$
As $q>p_i$, we have $\rho(n)\le\rho(n')$ in all cases. Putting
everything together, we see that
$$
\rho(2n_k)\log\log(2n_k)\le\rho(n)\log\log n
 \le\rho(n')\log\log n'.
$$
Since $n'\in\cS_{j+1}$ is arbitrary, we obtain~\eqref{eq:Sjbound}
with $j$ replaced by $j+1$, which completes the induction and
finishes our proof of~\eqref{eq:rho-ineq}.

Next, we turn to the proof of~\eqref{eq:lim-2nk}. Using the Prime
Number Theorem in the form
$$
\log n_k=\sum_{p\le p_k}\log p=(1+o(1))p_k\qquad(k\to\infty)
$$
together with Mertens' theorem
(see~\cite[Theorem~2.7(e)]{MontVau}), it is easy to see that
\begin{equation}
\label{eq:lim1} \lim_{k\to\infty}\left\{\log\log(2n_k)\prod_{p\le
p_k}\(1-\frac{1}{p}\)\right\}=e^{-\gamma}.
\end{equation}
Also,
\begin{equation}
\label{eq:lim2} \lim_{k\to\infty}\prod_{2<p\le
p_k}\frac{p(p-2)}{(p-1)^2}=\lim_{k\to\infty}
C_2\prod_{p>p_k}\(1+\frac{1}{p(p-2)}\)=C_2.
\end{equation}
By~\eqref{eq:rho-defn} we have
$$
\rho(2n_k)=\prod_{p\le p_k}\(1-\frac{1}{p}\)\prod_{2<p\le p_k}
\frac{p(p-2)}{(p-1)^2}\qquad(k\ge 1),
$$
and thus~\eqref{eq:lim-2nk} is an immediate consequence
of~\eqref{eq:lim1} and~\eqref{eq:lim2}.

\section{Proof of Theorem~\ref{thm:main}}
\label{sec:three}

As in~\cite[Th\'eor\`eme~3]{Nic} we put
$$
f(x)=e^\gamma\log\vartheta(x)\prod_{p\le
x}\(1-\frac{1}{p}\)\qquad(x\ge 2),
$$
where $\vartheta(x)=\sum_{p\le x}\log p$ is the Chebyshev
$\vartheta$-function. For our purposes, it is convenient to define
$$
g(x)=e^\gamma\log\(\vartheta(x)+\log 2\)\prod_{p\le x}
\(1-\frac{1}{p}\)\prod_{p>x}\(1+\frac{1}{p(p-2)}\)\qquad(x\ge 2),
$$
This definition is motivated by the fact that
$$
g(p_k)=C_2^{-1}e^\gamma\rho(2n_k)\log\log(2n_k)\qquad(k\ge 1).
$$
As mentioned earlier, the inequalities~\eqref{eq:RHtrue}
and~\eqref{eq:RHtrueB} are equivalent, and~\eqref{eq:RHtrueB} is
clearly equivalent to
$$
\log g(p_k)\le 0.
$$
Thus, to prove Theorem~\ref{thm:main} it suffices to study the
sign of $\log g(x)$.

By the trivial inequality $\log(1+t)\le t$ for all $t>-1$ and the
fact that $g(x)>f(x)$ for all $x\ge 2$, it is easy to see that
\begin{equation}
\label{eq:log g/f} 0< \log\frac{g(x)}{f(x)}\le\frac{\log
2}{\vartheta(x)\log\vartheta(x)}+\frac{1}{x-2}\qquad(x>2).
\end{equation}
Here, we have used the fact that
$$
\sum_{p>x}\frac{1}{p(p-2)}\le\sum_{n\ge\fl{x}+1}\frac{1}{n(n-2)}
=\frac{2\fl{x}-1}{2\fl{x}(\fl{x}-1)}<\frac{1}{x-2}\qquad(x>2).
$$

First, let us suppose that the Riemann hypothesis is true. In this
case, we have from~\cite[p.~383]{Nic}:
$$
\log f(x)\le -\frac{0.8}{\sqrt{x}\,\log x}\qquad(x\ge 3000).
$$
Using this bound in~\eqref{eq:log g/f} together with the
inequality $\vartheta(x)\ge 4x/5$ (which holds unconditionally for
$x\ge 121$ by~\cite[Theorems~4 and~18]{RosSch}), one sees that
$$
\log g(x)\le\frac{\log 2}{(4x/5)\log(4x/5)}+\frac{1}{x-2}-
\frac{0.8}{\sqrt{x}\,\log x}\le -\frac{0.6}{\sqrt{x}\,\log x}
$$
for all $x\ge 3000$.  This implies the desired
bound~\eqref{eq:RHtrueB} for all $k\ge 431$; for smaller values of
$k$, the bound~\eqref{eq:RHtrueB} may be verified by a direct
computation. This proves Theorem~\ref{thm:main} under the Riemann
hypothesis.

Next, suppose that the Riemann hypothesis is false, and let
$\theta$ denote the supremum of the real parts of the zeros of the
Riemann zeta function.  Then, by~\cite[Th\'eor\`eme~3c]{Nic} one
has
$$
\limsup_{x\to\infty} x^b\log f(x)>0 \qquad \text{and} \qquad
\liminf_{x\to\infty} x^b\log f(x)<0
$$
for any fixed number $b$ such that $1-\theta<b<1/2$. In view
of~\eqref{eq:log g/f} and the Chebyshev bound $\vartheta(x)\gg x$
it is clear that
$$
\log g(x)=\log f(x)+O(x^{-1});
$$
hence, we also have
$$
\limsup_{x\to\infty} x^b\log g(x)>0 \qquad \text{and} \qquad
\liminf_{x\to\infty} x^b\log g(x)<0.
$$
In particular, $\log g(p_k)$ changes sign infinitely often, which
implies Theorem~\ref{thm:main} if the Riemann hypothesis is
false.

\end{document}